\documentclass[twoside,leqno,10pt, A4]{amsart}
\usepackage{amsfonts}
\usepackage{amsmath}
\usepackage{amscd}
\usepackage{amssymb}
\usepackage{amsthm}
\usepackage{amsrefs}
\usepackage{latexsym}
\usepackage{mathrsfs}
\usepackage{bbm}
\usepackage{enumerate}
\usepackage{graphicx}

\usepackage{amsfonts}
\usepackage{amsmath}
\usepackage{amscd}
\usepackage{amssymb}
\usepackage{amsthm}
\usepackage{amsrefs}
\usepackage{latexsym}
\usepackage{mathrsfs}
\usepackage{bbm}
\usepackage{amscd}
\usepackage{amssymb}
\usepackage{amsthm}
\usepackage{amsrefs}
\usepackage{latexsym}
\usepackage{mathrsfs}
\usepackage{bbm}
\usepackage{enumerate}
\usepackage{graphicx}
\usepackage{color}
\begin{document}

\newtheorem{theorem}[subsection]{Theorem}
\newtheorem{proposition}[subsection]{Proposition}
\newtheorem{lemma}[subsection]{Lemma}
\newtheorem{corollary}[subsection]{Corollary}
\newtheorem{conjecture}[subsection]{Conjecture}
\newtheorem{prop}[subsection]{Proposition}
\numberwithin{equation}{section}
\newcommand{\mr}{\ensuremath{\mathbb R}}
\newcommand{\mc}{\ensuremath{\mathbb C}}
\newcommand{\dif}{\mathrm{d}}
\newcommand{\intz}{\mathbb{Z}}
\newcommand{\ratq}{\mathbb{Q}}
\newcommand{\natn}{\mathbb{N}}
\newcommand{\comc}{\mathbb{C}}
\newcommand{\rear}{\mathbb{R}}
\newcommand{\prip}{\mathbb{P}}
\newcommand{\uph}{\mathbb{H}}
\newcommand{\fief}{\mathbb{F}}
\newcommand{\majorarc}{\mathfrak{M}}
\newcommand{\minorarc}{\mathfrak{m}}
\newcommand{\sings}{\mathfrak{S}}
\newcommand{\fA}{\ensuremath{\mathfrak A}}
\newcommand{\mn}{\ensuremath{\mathbb N}}
\newcommand{\mq}{\ensuremath{\mathbb Q}}
\newcommand{\half}{\tfrac{1}{2}}
\newcommand{\f}{f\times \chi}
\newcommand{\summ}{\mathop{{\sum}^{\star}}}
\newcommand{\chiq}{\chi \bmod q}
\newcommand{\chidb}{\chi \bmod db}
\newcommand{\chid}{\chi \bmod d}
\newcommand{\sym}{\text{sym}^2}
\newcommand{\hhalf}{\tfrac{1}{2}}
\newcommand{\sumstar}{\sideset{}{^*}\sum}
\newcommand{\sumprime}{\sideset{}{'}\sum}
\newcommand{\sumprimeprime}{\sideset{}{''}\sum}
\newcommand{\sumflat}{\sideset{}{^\flat}\sum}
\newcommand{\shortmod}{\ensuremath{\negthickspace \negthickspace \negthickspace \pmod}}
\newcommand{\V}{V\left(\frac{nm}{q^2}\right)}
\newcommand{\sumi}{\mathop{{\sum}^{\dagger}}}
\newcommand{\mz}{\ensuremath{\mathbb Z}}
\newcommand{\leg}[2]{\left(\frac{#1}{#2}\right)}
\newcommand{\muK}{\mu_{\omega}}
\newcommand{\thalf}{\tfrac12}
\newcommand{\lp}{\left(}
\newcommand{\rp}{\right)}
\newcommand{\Lam}{\Lambda_{[i]}}
\newcommand{\lam}{\lambda}
\def\L{\fracwithdelims}
\def\om{\omega}
\def\pbar{\overline{\psi}}
\def\phi{\varphi}
\def\lam{\lambda}
\def\lbar{\overline{\lambda}}
\newcommand\Sum{\Cal S}
\def\Lam{\Lambda}
\newcommand{\sumtt}{\underset{(d,2)=1}{{\sum}^*}}
\newcommand{\sumt}{\underset{(d,2)=1}{\sum \nolimits^{*}} \widetilde w\left( \frac dX \right) }

\theoremstyle{plain}
\newtheorem{conj}{Conjecture}
\newtheorem{remark}[subsection]{Remark}

\providecommand{\re}{\mathop{\rm Re}}
\providecommand{\im}{\mathop{\rm Im}}
\def\cI{\mathcal{I}}
\def\cL{\mathcal{L}}
\def\E{\mathbb{E}}

\makeatletter
\def\widebreve{\mathpalette\wide@breve}
\def\wide@breve#1#2{\sbox\z@{$#1#2$}%
     \mathop{\vbox{\m@th\ialign{##\crcr
\kern0.08em\brevefill#1{0.8\wd\z@}\crcr\noalign{\nointerlineskip}%
                    $\hss#1#2\hss$\crcr}}}\limits}
\def\brevefill#1#2{$\m@th\sbox\tw@{$#1($}%
  \hss\resizebox{#2}{\wd\tw@}{\rotatebox[origin=c]{90}{\upshape(}}\hss$}
\makeatletter

\title[Lower bounds for moments of the derivative of the Riemann zeta function]{Lower bounds for moments of the derivative of the Riemann zeta function}

\author{Peng Gao}
\address{School of Mathematical Sciences, Beihang University, Beijing 100191, P. R. China}
\email{penggao@buaa.edu.cn}
\begin{abstract}
 We establish in this paper sharp lower bounds for the $2k$-th moment of the derivative of the Riemann zeta function on the critical line for all real $k \geq 0$.
\end{abstract}

\maketitle

\noindent {\bf Mathematics Subject Classification (2010)}: 11M06 \newline

\noindent {\bf Keywords}: lower bounds, moments, derivative,  Riemman zeta function

\section{Introduction}
\label{sec 1}

  It is an important subject in analytical number theory to investigate moments of the Riemann zeta function $\zeta(s)$ on the critical line as they can be applied to study the maximum size of the zeta function as well as primes in short intervals via zero density estimates. We denote the $2k$-th moment of $\zeta(s)$ on the critical line by
\begin{align*}
 M_{k}(T)=\int\limits^{2T}_{T}| \zeta(\half+it)|^{2k} \dif t.
\end{align*}
  The study on $M_k(T)$ dates back to the work of G.H. Hardy and J. E. Littlewood \cite{H&L}, who established an asymptotic formula for $M_1(T)$. In \cite{In},  A. E. Ingham established an asymptotic formula for $M_2(T)$. No other asymptotic formulas are known for $M_k(T)$ except for the trivial case $k=0$. Despite of this, J. P. Keating and N. C. Snaith \cite{Keating-Snaith02} made precise conjectured formulas for $M_k(T)$ for all real $k \geq 0$ by drawing analogues with the random matrix theory.
Using multiple Dirichlet series, A. Diaconu, D. Goldfeld and J. Hoffstein \cite{DGH} also obtained the same conjectured formulas. More precise asymptotic formulas with lower order terms were conjectured by J. B. Conrey, D. W. Farmer, J. P. Keating, M. O. Rubinstein and N. C. Snaith in \cite{CFKRS}.

  Owing much to the work in \cites{Ramachandra2, Ramachandra1, Ramachandra3, H-B81-2, Radziwill, BCR, R&Sound, Radziwill&Sound, Sound2009, Harper, HRS, H&Sound, CG, Sound95, BR, CG92, CG01}, we now have sharp upper and lower bounds for $M_k(T)$ of the conjectured order of magnitude for all $k \geq 0$ with some of them being valid under the truth of the Riemann hypothesis (RH).

  Among the many methods applied in the above work, we point out notably a simple and powerful method developed by Z. Rudnick and K. Soundararajan  \cite{R&Sound, R&Sound1} towards establishing sharp lower bounds for moments of families of $L$-functions, a method of K. Soundararajan \cite{Sound2009} and its refinement by A. J. Harper \cite{Harper} to derive sharp upper bounds for moments of families of $L$-functions under the generalized Riemann hypothesis (GRH). We note also an upper bounds principle developed by M. Radziwi{\l\l} and K. Soundararajan in \cite{Radziwill&Sound} for establishing upper bounds for moments of families of $L$-functions as well as its dual lower bounds principle developed by W. Heap and K. Soundararajan in \cite{H&Sound}.

  Similar to $M_k(T)$, it is also interesting to study moments of the derivatives of $\zeta(s)$ on the critical line. For integers $l \geq 1$, let
\begin{align*}
 I_{k,l}(T)=\int\limits^{T}_{1}| \zeta^{(l)}(\half+it)|^{2k} \dif t.
\end{align*}

  An asymptotic formula for $I_{1, l}(T)$ is also given in the above mentioned work of A. E. Ingham \cite{In}. In \cite{Conrey1988}, J. B. Conrey  obtained an asymptotic formula for $I_{2,l}(T)$. Also, in connection with the random matric theory,  J.B. Conrey, M.O. Rubinstein and N.C. Snaith \cite[Conjecture 1]{CRS06} conjectured that
\begin{align*}
 I_{k,1}(T) \sim a_kb_kT(\log T)^{k^2+2k},
\end{align*}
   for some explicit constants $a_k, b_k$.

   Under RH, M. B. Milinovich \cite{Milinovich2011} established essentially upper bounds of the correct order of magnitude for $I_{k,l}(T)$ for positive integers $k, l$. His result was further improved to yield optimal upper bounds by A. Ivi\'{c} \cite{Ivic17} for $I_{k,2}(T)$ for positive integers $k$. The methods employed in \cite{Milinovich2011} and \cite{Ivic17} allow one to deduce upper bounds for $I_{k,l}(T)$ from the corresponding ones for $M_{k}(T)$. As sharp upper bounds for $M_k(T)$ are known for all $k \geq 0$ under RH from the work of K. Soundararajan \cite{Sound2009} and for all $0 \leq k \leq 2$ unconditionally from the work of W. Heap, M. Radziwi{\l\l} and K. Soundararajan \cite{HRS}, we may apply the methods in  \cite{Milinovich2011} and \cite{Ivic17} to derive that unconditionally for $1/2 \leq k \leq 2$ and under RH for $k \geq 2$, we have for all integers $l \geq 1$,
\begin{align*}
  I_{k,l}(T)\ll_{k, l} T(\log T)^{k^2+2kl}.
\end{align*}

  On the other hand, K. Sono \cite{Sono12}, T. Christ and J. Kalpokas \cite{CK} studied lower bounds for $I_{k,l}(T)$. It follows from \cite[Corollary 1.1]{CK} that we have for any rational $k \geq 1$ and any positive integer $l$,
\begin{align*}
  I_{k,l}(T)\gg_{k, l} T(\log T)^{k^2+2kl}.
\end{align*}

 The aim of this paper is to obtain sharp lower bounds for $I_{k,l}(T)$ for all real $k \geq 0$. For simplicity, we shall focus on $I_{k, 1}(T)$ throughout, although our methods carry over to treat $I_{k, l}(T)$ for other $l$ as well.  Our main result is as follows.
\begin{theorem}
\label{thmlowerboundJ}
   For large $T$ and any $k \geq 0$, we have
\begin{align*}
   I_{k,1}(T) \gg_{k} T(\log T)^{k^2+2k}.
\end{align*}
\end{theorem}

  Combining Theorem \ref{thmlowerboundJ} and our discussions above, we immediately obtain the following result concerning the order of magnitude of $I_{k,1}(T) $.
\begin{corollary}
\label{cororderofmag}
   For large $T$ and any $1/2 \leq k \leq 2$, we have
\begin{align*}
   I_{k,1}(T) \asymp T(\log T)^{k^2+2k}.
\end{align*}
\end{corollary}

  The proof of Theorem \ref{thmlowerboundJ} is based on the above mentioned lower bounds principle of W. Heap and K. Soundararajan \cite{H&Sound}, together with the approach taken in \cite[Section 5]{CK}. We note that a similar approach to the one used in \cite[Section 5]{CK} has already been employed by
M. B. Milinovich and N. Ng \cite{MN} in their study on lower bounds for the discrete moments of the derivative of
$\zeta(s)$ at nontrivial zeros. Since these discrete moments can be regarded as analogues to $I_{k,1}(T)$ and are studied by the author in \cite{Gao2021-5}, the proof of Theorem \ref{thmlowerboundJ} also makes use of some approaches there as well.

\section{Preliminaries}
\label{sec 2}

 We reserve the letter $p$ for a prime number in this paper and we recall the following well-known results on sums of primes (see \cite[Theorem 2.7]{MVa1}).
\begin{lemma} \label{RS} Let $x \geq 2$. We have, for some constant $b$,
\begin{align*}
&\sum_{p\le x} \frac{1}{p} = \log \log x + b+ O\Big(\frac{1}{\log x}\Big), \\
&\sum_{p\le x} \frac {\log p}{p} = \log x + O(1).
\end{align*}
\end{lemma}

  We note the following mean value theorem given in \cite[Lemma 4.1]{MN1} concerning integrals over Dirichlet polynomials.
\begin{lemma}
\label{Lem-MVDP}
	Let $\{a_n\}$ and $\{b_n\}$ be sequences of complex numbers. Let $T_1$ and $T_2$ be positive real numbers and $g(t)$ be a real-valued function that is continuously differentiable on the interval $[T_1, T_2]$. Then
\begin{align*}
\begin{split}
& \int^{T_2}_{T_1}g(t)\left ( \sum^{\infty}_{n=1}a_nn^{-it}\right )\left( \sum^{\infty}_{n=1}b_nn^{it}\right )dt \\
=& \int^{T_2}_{T_1}g(t)dt\sum^{\infty}_{n=1}a_nb_n+O\left ( \left (|g(T_1)|+|g(T_2)|+\int^{T_2}_{T_1}|g'(t)|dt \right )\left ( \sum^{\infty}_{n=1}n|a_n|^2\right )^{1/2}\left( \sum^{\infty}_{n=1}n|b_n|^2 \right )^{1/2}\right ) .
\end{split}
\end{align*}
\end{lemma}

\section{Proof of Theorem \ref{thmlowerboundJ}}
\label{sec 2'}

\subsection{The lower bound principle}

    We may assume that $k>0$ and $T$ is a large number throughout the proof.
We also point out that the explicit constants involved in estimations using $\ll$ or the big-$O$ notations in the proof depend on $k$ only and are uniform with
 respect to $p$ and $T$.

 We follow the ideas of A. J. Harper in \cite{Harper} to define for a large number $M$ depending on $k$ only,
$$ \alpha_{0} = 0, \;\;\;\;\; \alpha_{j} = \frac{20^{j-1}}{(\log\log T)^{2}} \;\;\; \forall \; j \geq 1, \quad
\mathcal{J} = \mathcal{J}_{k,T} =\max\{j : \alpha_{j} \leq 10^{-M} \} . $$

   We denote $\ell_j :=\lceil e^2k\alpha^{-3/4}_j \rceil$ for $1 \leq j \leq \mathcal{J}$ and divide the interval $(0, T^{\alpha_{\mathcal{J}}}]$ into disjoint subintervals $I_j=(T^{\alpha_{j-1}}, T^{\alpha_{j}}], 1 \leq j \leq \mathcal{J}$. We define for any real number $\ell$ and any $x \in \mr$,
\begin{align*}
  E_{\ell}(x) = \sum_{j=0}^{\lceil \ell \rceil} \frac {x^{j}}{j!}.
\end{align*}

 We also define for any real number $\alpha$ and any $1\leq j \leq \mathcal{J}$,
\begin{align*}
 {\mathcal P}_j(s)=&  \sum_{ p \in I_j}  \frac{1}{p^s}, \quad {\mathcal N}_j(s, \alpha) = E_{\ell_j } \Big (\alpha {\mathcal P}_j(s) \Big ), \quad  {\mathcal N}(s, \alpha)=  \prod^{\mathcal{J}}_{j=1} {\mathcal N}_j(s, \alpha).
\end{align*}

 We deduce from \cite[(3.1)]{Gao2021-5} and Lemma \ref{RS} that for any large number $N$, we can take $T$ large enough so that
\begin{align}
\label{sumpj}
 {\mathcal P}_1(1) \leq  \frac 1N \ell_1, \quad {\mathcal P}_j(1) \leq  \min (10, \frac 1N \ell_j), \quad 2 \leq j \leq \mathcal{J}.
\end{align}

 We denote $\Omega(n)$ for the number of prime powers dividing $n$ and $g(n)$ for the multiplicative function given on prime powers by $g(p^{r}) = 1/r!$ and define functions $b_j(n), 1 \leq j \leq {\mathcal{J}}$  such that $b_j(n)=0$ or $1$ and that $b_j(n)=1$ only when $\Omega(n) \leq \ell_j$ and all the prime factors of $n$ are from the interval $I_j$. We then have
\begin{align*}
 {\mathcal N}_j(s,\alpha) = \sum_{n_j}  \frac{\alpha^{\Omega(n_j)}}{g(n_j)}  b_j(n_j) \frac 1{n^s_j}, \quad 1\le j \le {\mathcal{J}}.
\end{align*}

  It follows from \cite[Section 3.1]{Gao2021-5} that each ${\mathcal N}_j(s,\alpha)$ is a short Dirichlet polynomial of length at most $T^{\alpha_{j}\lceil e^2k\alpha^{-3/4}_j \rceil}$ and that ${\mathcal N}(s, \alpha)$ is also a short Dirichlet polynomial of length at most $T^{40 e^2k10^{-M/4}}$.

  We now write for simplicity that
\begin{align}
\label{Nskexpression}
 {\mathcal N}(s, \alpha)= \sum_{n} \frac{a_{\alpha}(n)}{n^s}.
\end{align}
   We note that $a_{\alpha}(n) \neq 0$ only when $n=\prod_{1\leq j \leq \mathcal{J}}n_j$ such that $b_j(n_j)=1$, in which case we have
\begin{align*}
  a_{\alpha}(n)= \prod_{n_j}  \frac{\alpha^{\Omega(n_j)}}{g(n_j)}  b_j(n_j).
\end{align*}

   We combine \cite[(2.1), (3.3)]{Gao2021-5} to see that for all $n \geq 3$,
\begin{align}
\label{anbound}
  a_{\alpha}(n) \leq e^{\frac {|\alpha|\log n}{\log \log n}(1+O(\frac {1}{\log \log n}))} \text{  and} \quad a_k(n) = 0 \text{ when } n > T^{40 e^2k10^{-M/4}}.
\end{align}

   Moreover, we note that \cite[(3.4)]{Gao2021-5} implies that for $\Re(s) \geq -1/\log T$ and $T$ large enough,
\begin{align}
\label{Nnormbound}
 |{\mathcal N}(s,\alpha)| \ll e^{|\alpha|\frac {\log T}{\log \log T}(1+O(\frac {1}{\log \log T}))}T^{40e^2k10^{-M/4}(1+1/\log T)}.
\end{align}

  In the proof of Theorem \ref{thmlowerboundJ}, we need the following lower bounds principle of W. Heap and K. Soundararajan
  \cite{H&Sound} for our case.
\begin{lemma}
\label{lem1}
 With notations as above, we have for $0<k \leq 1/2$,
\begin{align}
\label{basiclowerbound}
\begin{split}
\int^T_1-\zeta'(\half+it)\mathcal{N}(\half+it, k-1)\mathcal{N}(\half-it, k)dt
 \ll & \Big ( \int^T_1|\zeta'(\half+it)|^{2k}dt \Big )^{1/2}\Big ( \int^T_1|\zeta'(\half+it)|^{2} |\mathcal{N}(\half+it, k-1)|^2 dt \Big)^{(1-k)/2} \\
 & \times \Big ( \int^T_1  \prod^{\mathcal{J}}_{j=1}\big ( |{\mathcal N}_j(\half+it, k)|^2+ |{\mathcal Q}_j(\half+it,k)|^{2r_k} \big )dt
 \Big)^{k/2}.
\end{split}
\end{align}
 Also, we have for $k>1/2$,
\begin{align}
\label{basicboundkbig}
\begin{split}
 & \int^T_1-\zeta'(\half+it)\mathcal{N}(\half+it, k-1)\mathcal{N}(\half-it, k)dt \\
 \leq & \Big (  \int^T_1|\zeta'(\half+it)|^{2k}dt \Big )^{\frac {1}{2k}}\Big (  \int^T_1\prod^{\mathcal{J}}_{j=1} \big ( |{\mathcal N}_j(\half+it, k)|^2+ |{\mathcal Q}_j(\half+it,k)|^{2r_k} \big )dt \Big)^{\frac {2k-1}{2k}}.
\end{split}
\end{align}
  Here the implied constants in \eqref{basiclowerbound} and \eqref{basicboundkbig} depend on $k$ only,  and we define
$$
{\mathcal Q}_j(s,k) =\Big( \frac{64 \max (2, k+3/2 ) {\mathcal P}_j(s)}{\ell_j} \Big)^{\ell_j},
$$
  with $r_k=2+\lceil 1/k \rceil$ for $0<k \leq 1/2$ and $r_k=1+\lceil 2k/(2k-1) \rceil$ for $k >1/2$.
\end{lemma}

  We skip the proof of the above lemma as it can be established similar to those of \cite[Lemma 3.2-3.3]{Gao2021-5}. We deduce from the above lemma that in order to establish Theorem \ref{thmlowerboundJ}, it suffices to prove the following three propositions.
\begin{proposition}
\label{Prop4}
  With notations as above, we have for $k > 0$,
\begin{align}
\label{L1estmation}
 \int^T_1-\zeta'(\half+it)\mathcal{N}(\half+it, k-1)\mathcal{N}(\half-it, k)dt \gg T (\log T)^{ k^2+1}.
\end{align}
\end{proposition}

\begin{proposition}
\label{Prop6}
  With notations as above, we have for $0<k \leq 1/2$,
\begin{align}
\label{L2estmation}
\int^T_1|\zeta'(\half+it)|^{2} |\mathcal{N}(\half+it, k-1)|^2 dt   \ll T ( \log T  )^{k^2+2}.
\end{align}
\end{proposition}

\begin{proposition}
\label{Prop5}
  With notations as above, we have for $k >0$,
\begin{align*}
\int^T_1  \prod^{\mathcal{J}}_{i=1}\big ( |{\mathcal N}_i(\half+it, k)|^2+ |{\mathcal Q}_i(\half+it,k)|^{2r_k} \big )dt \ll & T ( \log T  )^{k^2}.
\end{align*}
\end{proposition}

  We shall omit the proof of Proposition \ref{Prop5} as it is similar to that of \cite[Proposition 3.5]{Gao2021-5}, upon making use of Lemma \ref{Lem-MVDP}. In the rest of the paper, we shall prove the remaining two propositions.

\subsection{Proof of Proposition \ref{Prop4}}
\label{sec 4.8}

   The proof is based on the approaches used in Section 5 of \cite{CK} and Section 5 of \cite{MN1}. We denote the left side expression in \eqref{L1estmation} by $S_1$ and apply Cauchy’s residue theorem to deduce that
\begin{align*}
   S_1 &  = \frac{1}{2 \pi i} \int_{\mathcal{C}} -\zeta'(s)  \mathcal{N}(s, k-1)\mathcal{N}(1-s, k) \, ds,
\end{align*}
where $\mathcal{C}$ consists of line segments from $\half+i$ to $\kappa+i$, then from $ \kappa+i$ to $\kappa+iT$ and lastly from $\kappa+iT$ to $\half+iT$, where $\kappa=1+(\log T)^{-1}$.

 We apply \eqref{Nnormbound} and the estimation (see \cite[(20)]{Gonek})
\begin{equation*}
\begin{split}
   \zeta'(s)  \ll
\left\{
 \begin{array}
  [c]{ll}
   (1+|t|)^{(1-\Re(s))/2+\epsilon}, \quad 0 \leq \Re(s) \leq 1, \\
   (1+|t|)^{\epsilon} , \quad \Re(s) \geq 1,
\end{array}
\right .
\end{split}
\end{equation*}
  to see that the integral is bounded by $O(T^{1-\varepsilon})$ on the horizontal edges of the contour. We thus conclude that
\begin{align}
\label{Sexp}
   S_1 = S_{1,R}+ O(T^{1-\varepsilon}),
\end{align}
  where
\begin{align}
\label{SRint}
  S_{1,R} = \frac{1}{2 \pi i} \int_{\kappa+i}^{\kappa+iT} -\zeta'(s)
 \mathcal{N}(s, k-1)\mathcal{N}(1-s, k) \, ds+O(T^{1-\varepsilon}).
\end{align}

  To evaluate $S_{1, R}$, we define the Dirichlet convolution $f*g$ for two arithmetic functions $f(k), g(k)$ by
\begin{align*}
  f*g(k)=\sum_{mn=k}f(m)g(n).
\end{align*}

   Using this notation and that given in \eqref{Nskexpression}, we apply Lemma \ref{Lem-MVDP} to evaluate $S_{1,R}$ in \eqref{SRint} to see that
\begin{align*}
  S_{1,R} =&  \frac{T-1}{2 \pi} \sum_{n}\frac {(\log * a_{k-1})(n) \cdot a_k(n)}{n}+ O\left (  \left (\sum^{\infty}_{n=1} \frac {(\log * a_{k-1})(n)^2}{n^{2\kappa-1}}\right )^{\half}\left (\sum^{\infty}_{n=1} \frac {a_{k}(n)^2}{n^{1-2\kappa}}\right )^{\half}\right ).
\end{align*}

  We apply the estimations given in \eqref{anbound} to see that for $T$ large enough,
\begin{align*}
 \sum^{\infty}_{n=1} \frac {a_{k}(n)^2}{n^{1-2\kappa}} \ll e^{4k\log T/\log \log T}\sum_{n \leq T^{40 e^2k10^{-M/4}}} \frac {1}{n^{1-2\kappa}}\ll T^{1-\varepsilon}.
\end{align*}

  Moreover, we have that
\begin{align*}
\begin{split}
 & (\log * a_{k-1}) (n) \leq \log n \sum_{n \leq T^{40 e^2k10^{-M/4}}}|a_{k-1}(n)| \leq T^{1/2-\varepsilon}\log n.
\end{split}
\end{align*}
  It follows from the above that
\begin{align*}
  \sum^{\infty}_{n=1} \frac {(\log * a_{k-1})(n)^2}{n^{2\kappa-1}} \ll T^{1-2\varepsilon}\sum^{\infty}_{n=1} \frac {\log^2 n}{n^{2\kappa-1}} \ll T^{1-\varepsilon},
\end{align*}
   where the last estimation above follows from the bound that (see \cite[(16)]{MN1}) uniformly for $\sigma>1$ and any integer $i \geq 0$,
\begin{align*}
  \sum^{\infty}_{n=1} \frac {\log^i n}{n^{\sigma}} \ll \frac 1{(\sigma-1)^{i+1}}.
\end{align*}

  We then conclude from the above discussions that
\begin{align}
\label{SRexp}
\begin{split}
  S_{1, R} =&  \frac{T-1}{2\pi }\sum_{n, m}\frac {a_{k-1}(m)a_k(mn)(\log n)}{mn}+O(T^{1-\varepsilon}) \\
=& \frac{T-1}{2\pi }\sum_{n}\frac {\log n}{n}\sum_{m}\frac {a_{k-1}(m)a_k(mn)}{m}+O(T^{1-\varepsilon}).
\end{split}
\end{align}

  It remains to estimate the last expression above. To do so, we may assume that $n=\prod^{\mathcal{J}}_{j=1}n_j$ with $b_j(n_j)=1$ for $1 \leq j \leq \mathcal{J}$.  Then the inner sum of the last expression above becomes
\begin{align*}
\begin{split}
 \sum_{m}\frac {a_{k-1}(m)a_k(mn)}{m}=& \prod^{\mathcal{J}}_{j=1}\Big ( \sum_{m_j} \frac{1}{m_j} \frac{k^{\Omega(n_jm_j)}(k-1)^{\Omega(m_j)}}{g(n_jm_j)g(m_j)}
b_j(n_jm_j)b_j(m_j) \Big ) \\
=& \prod^{\mathcal{J}}_{j=1}\Big ( \sum_{m_j} \frac{1}{m_j} \frac{k^{\Omega(n_jm_j)}(k-1)^{\Omega(m_j)}}{g(n_jm_j)g(m_j)}
b_j(n_jm_j) \Big ),
\end{split}
\end{align*}
  where the last equality above follows by noting that $b_j(n_jm_j)=1$ implies that $b_j(m_j)=1$ for all $1 \leq j \leq \mathcal{J}$.

Note that the factor $b_j(n_jm_j)$ restricts $m_j$ to have all prime factors in $I_j$ such that $\Omega(n_jm_j) \leq \ell_j$. If we remove this restrictions on $\Omega$, then the sum over $m_j$ becomes
\begin{align*}
\begin{split}
\sum_{m_j} \frac{1}{m_j} \frac{k^{\Omega(n_jm_j)}(k-1)^{\Omega(m_j)}}{g(n_jm_j)g(m_j)}
=\prod_{\substack{p\in I_j \\ (p, n_j)=1}}\Big (1+ \frac {k(k-1)}p+O(\frac 1{p^2}) \Big )\prod_{\substack{p_{i,j} \in I_j \\ p^{l_{i,j}}_{i,j} \| n_j\\ l_{i,j} \geq 1}}\Big (\frac {k^{l_{i,j}}}{l_{i,j}!}+ \frac {k^{l_{i,j}+1}(k-1)}{(l_{i,j}+1)!p_{i,j}}+ \frac {k^{l_{i,j}+2}(k-1)^2}{(l_{i,j}+2)!2!p^2_{i,j}}+\cdots  \Big ).
\end{split}
\end{align*}
  We recast the last product above as
\begin{align*}
\begin{split}
\prod_{\substack{p_{i,j} \in I_j \\ p^{l_{i,j}}_{i,j} \| n_j \\ l_{i,j} \geq 1}}\Big (\frac {k^{l_{i,j}}}{l_{i,j}!}+ \frac {k^{l_{i,j}+1}(k-1)}{(l_{i,j}+1)!p_{i,j}}+ \frac {k^{l_{i,j}+2}(k-1)^2}{(l_{i,j}+2)!2!p^2_{i,j}}+\cdots  \Big )
= \frac {k^{\Omega(n_j)}}{g(n_j)}\prod_{\substack{p \in I_j \\ p | n_j}}\Big (1 +O(\frac 1p) \Big ),
\end{split}
\end{align*}
   and we note that each factor in the last product above is positive.

  Using Rankin's trick by noticing that $2^{\Omega(n_jm_j)-\ell_j}\ge 1$ if $\Omega(n_jm_j) > \ell_j$,  we see that the error
   introduced this way does not exceed
\begin{align*}
\begin{split}
 & \sum_{m_j} \frac{1}{m_j} \frac{k^{\Omega(n_jm_j)}|1-k|^{\Omega(m_j)}}{g(n_jm_j)g(m_j)}2^{\Omega(n_jm_j)-\ell_j}
  \\
= & 2^{\Omega(n_j)-\ell_j} \sum_{m_j} \frac{1}{m_j} \frac{k^{\Omega(n_jm_j)}2^{\Omega(m_j)}|1-k|^{\Omega(m_j)}}{g(n_jm_j)g(m_j)} \\
\leq & 2^{\Omega(n_j)-\ell_j/2}\frac {k^{\Omega(n_j)}}{g(n_j)} \prod_{\substack{p\in I_j \\ (p, n_j)=1}}\Big (1+ \frac {k(k-1)}p+O(\frac 1{p^2}) \Big )\prod_{\substack{p \in I_j \\ p | n_j}}\Big (1 +O(\frac 1p) \Big ),
\end{split}
\end{align*}
  where the last estimation above follows from \eqref{sumpj}.

  It follows that we may write
\begin{align*}
\begin{split}
  \sum_{m}\frac {a_{k-1}(m)a_k(mn)}{m} = & \prod_{\substack{p\in \bigcup^{\mathcal{J}}_{j=1} I_j}}\Big (1+ \frac {k(k-1)}p+O(\frac 1{p^2}) \Big ) \\
& \times \prod^{\mathcal{J}}_{j=1}\Big(1+ f_j(n_j) \Big )\frac {k^{\Omega(n_j)}}{g(n_j)}\prod_{\substack{p \in I_j \\ p | n_j}}\Big (1 +O(\frac 1p) \Big )\Big (1+ \frac {k(k-1)}{p}+O(\frac 1{p^2}) \Big )^{-1} \\
= & \prod_{\substack{p\in \bigcup^{\mathcal{J}}_{j=1} I_j}}\Big (1+ \frac {k(k-1)}p+O(\frac 1{p^2}) \Big ) \prod^{\mathcal{J}}_{j=1}\Big(1+ f_j(n_j)\Big ) \frac {k^{\Omega(n_j)}}{g(n_j)} \prod_{\substack{p \in I_j \\ p | n_j}}\Big (1 +O(\frac 1p) \Big ),
\end{split}
\end{align*}
  where
\begin{align*}
\begin{split}
  | f_j(n_j)| \leq 2^{\Omega(n_j)-\ell_j/2}.
\end{split}
\end{align*}

  We apply the above estimation to see that
\begin{align}
\label{loglowerbound}
\begin{split}
 &  \sum_{n}\frac {\log n}{n}\sum_{m}\frac {a_{k-1}(m)a_k(mn)}{m} \\
= & \prod_{\substack{p\in \bigcup^{\mathcal{J}}_{j=1} I_j}}\Big (1+ \frac {k(k-1)}p+O(\frac 1{p^2}) \Big )\sum_{\substack{n=\prod_jn_j  }}\frac {\log n}{n}\prod^{\mathcal{J}}_{j=1}\Big(1+ f_j(n_j)\Big ) \frac {k^{\Omega(n_j)}}{g(n_j)}b_j(n_j)   \prod_{\substack{p \in I_j \\ p | n_j}}\Big (1 +O(\frac 1p) \Big ).
\end{split}
\end{align}

  Note that we have
\begin{align*}
\begin{split}
 & \sum_{\substack{n=\prod_jn_j}}\frac {\log n}{n}\prod^{\mathcal{J}}_{j=1}\Big(1+ f_j(n_j)\Big ) \frac {k^{\Omega(n_j)}}{g(n_j)} b_j(n_j)  \prod_{\substack{p \in I_j \\ p | n_j}}\Big (1 +O(\frac 1p) \Big ) \\
=& \sum_{\substack{n=\prod_jn_j}}\prod^{\mathcal{J}}_{j=1}\frac 1{n_j}\Big(1+ f_j(n_j) \Big )\frac {k^{\Omega(n_j)}}{g(n_j)} b_j(n_j) \prod_{\substack{p \in I_j \\ p | n_j}}\Big (1 +O(\frac 1p) \Big )\Big(\sum_{j}\log n_{j} \Big ) \\
=& \sum^{\mathcal{J}}_{j'=1} \prod^{\mathcal{J}}_{\substack{j=1\\ j \neq j'}}\Big ( \sum_{\substack{n_j}}\frac 1{n_j}\Big(1+ f_j(n_j) \Big ) \frac {k^{\Omega(n_j)}}{g(n_j)}b_j(n_j) \prod_{\substack{p \in I_j \\ p | n_j}}\Big (1 +O(\frac 1p) \Big ) \Big )\\
& \times \Big (  \sum_{\substack{ n_{j'}}}\frac {\log n_{j'}}{n_{j'}}\Big(1+ f_{j'}(n_{j'}) \Big ) \frac {k^{\Omega(n_{j'})}}{g(n_{j'})} b_{j'}(n_{j'}) \prod_{\substack{p \in I_{j'} \\ p | n_{j'}}}\Big (1 +O(\frac 1p) \Big ) \Big ).
\end{split}
\end{align*}

   We denote $\mathcal{N}_j, 1 \leq j \leq \mathcal{J}$ for the set of integers $n_j$ such that $n_j$ is divisible only by primes $p \in I_j$. We estimate the last sum of the last expression above by observing that $1-2^{\Omega(n_{j})-\ell_{j}/2} \leq 0$ when $\Omega(n_{j})\geq  \ell_{j}/2$, so that
\begin{align*}
& \sum_{\substack{ n_j }}\frac {\log n_{j}}{n_{j}}\Big(1+ f_j(n_j) \Big ) \frac {k^{\Omega(n_{j})}}{g(n_{j})} b_j(n_j) \prod_{\substack{p \in I_j \\ p | n_j}}\Big (1 +O(\frac 1p) \Big ) \geq  \sum_{\substack{n_j \in \mathcal{N}_j }}\frac {\log n_{j}}{n_{j}}\Big(1-2^{\Omega(n_{j})-\ell_{j}/2} \Big ) \frac {k^{\Omega(n_{j})}}{g(n_{j})} \prod_{\substack{p \in I_j \\ p | n_j}}\Big (1 +O(\frac 1p) \Big ).
\end{align*}

  We further observe that
\begin{align*}
 & \sum_{\substack{n_j \in \mathcal{N}_j }}\frac {\log n_{j}}{n_{j}}\Big(1-  2^{\Omega(n_{j})-\ell_{j}/2} \Big ) \frac {k^{\Omega(n_{j})}}{g(n_{j})} \prod_{\substack{p \in I_j \\ p | n_j}}\Big (1 +O(\frac 1p) \Big )\\
=& -\frac {d}{d s}\Big ( \sum_{\substack{n_j \in \mathcal{N}_j }}\frac {1}{n^{1+s}_{j}}\Big(1- 2^{\Omega(n_{j})-\ell_{j}/2} \Big )\frac {k^{\Omega(n_{j})}}{g(n_{j})} \prod_{\substack{p \in I_j \\ p | n_j}}\Big (1 +O(\frac 1p) \Big ) \Big ) \Big |_{s=0}.
\end{align*}

    Upon writing
\begin{align*}
  \sum_{\substack{n_j \in \mathcal{N}_j }}\frac {1}{n^{1+s}_{j}}\frac {k^{\Omega(n_{j})}}{g(n_{j})} \prod_{\substack{p \in I_j \\ p | n_j}}\Big (1 +O(\frac 1p) \Big ) =\prod_{\substack{p \in I_{j}}}\Big (1+\frac {k}{p^{1+s}}\Big (1+ O(\frac 1{p}) \Big )+\frac {k^2}{2!p^{2(1+s)}}\Big (1+ O(\frac 1{p}) \Big )+\cdots \Big ),
\end{align*}
   we deduce that
\begin{align*}
 & -\frac {d}{d s}\Big ( \sum_{\substack{n_j \in \mathcal{N}_j }}\frac {1}{n^{1+s}_{j}}\frac {k^{\Omega(n_{j})}}{g(n_{j})} \prod_{\substack{p \in I_j \\ p | n_j}}\Big (1 +O(\frac 1p) \Big ) \Big ) \Big |_{s=0} \\
=& \prod_{\substack{p \in I_{j}}}\Big (1+\frac {k}{p}+ O(\frac 1{p^2}) \Big ) \Big (\sum_{\substack{p \in I_{j}}}\Big (\frac {k\log p}{p}+O(\frac 1{p^2})\Big )\Big (1+O(\frac 1{p}) \Big )^{-1}\Big ) \\
=& \prod_{\substack{p \in I_{j}}}\Big (1+\frac {k}{p}+ O(\frac 1{p^2}) \Big ) \Big (\sum_{\substack{p \in I_{j}}}\frac {k\log p}{p}+O(\frac 1{p^2})\Big ) .
\end{align*}

  Note also that we have
\begin{align*}
 & -\frac {d}{d s}\Big (\sum_{\substack{n_j \in \mathcal{N}_j }}\frac {2^{\Omega(n_{j})-\ell_{j}/2} }{n^{1+s}_{j}}\frac {k^{\Omega(n_{j})}}{g(n_{j})} \prod_{\substack{p \in I_j \\ p | n_j}}\Big (1 +O(\frac 1p) \Big ) \Big ) \Big |_{s=0} \\
=& 2^{-\ell_{j}/2}\prod_{\substack{p \in I_{j}}}\Big (1+\frac {2k}{p}+ O(\frac 1{p^2}) \Big ) \Big (\sum_{\substack{p \in I_{j}}}\frac {2k\log p}{p}+O(\frac 1{p^2})\Big ) \\
\leq &  2^{-\ell_{j}/4}\prod_{\substack{p \in I_{j}}}\Big (1+\frac {k}{p}+ O(\frac 1{p^2}) \Big ) \Big (\sum_{\substack{p \in I_{j}}}\frac {k\log p}{p}+O(\frac 1{p^2}) \Big ).
\end{align*}

  It follows that
\begin{align*}
& \sum_{\substack{n_j \in \mathcal{N}_j }}\frac {\log n_{j}}{n_{j}}\Big(1- 2^{\Omega(n_{j})-\ell_{j}/2} \Big ) \frac {k^{\Omega(n_{j})}}{g(n_{j})} \prod_{\substack{p \in I_j \\ p | n_j}}\Big (1 +O(\frac 1p) \Big ) \\
\geq &  (1-2^{-\ell_{j}/4})\prod_{\substack{p \in I_{j}}}\Big (1+\frac {k}{p}+ O(\frac 1{p^2}) \Big ) \Big (\sum_{\substack{p \in I_{j}}}\frac {k\log p}{p}+O(\frac 1{p^2})\Big ).
\end{align*}

   We apply similar arguments as above to see that we have
\begin{align*}
& \sum_{\substack{n_j}}\frac 1{n_j}\Big(1+ f_j(n_j) \Big ) \frac {k^{\Omega(n_j)}}{g(n_j)}b_j(n_j) \prod_{\substack{p \in I_j \\ p | n_j}}\Big (1 +O(\frac 1p) \Big ) \\
\geq & \sum_{\substack{n_j \in \mathcal{N}_j }}\frac {1}{n_{j}}\Big(1- 2^{\Omega(n_{j})-\ell_{j}/2}\Big ) \frac {k^{\Omega(n_{j})}}{g(n_{j})} \prod_{\substack{p \in I_j \\ p | n_j}}\Big (1 +O(\frac 1p) \Big ) \geq   (1-2^{-\ell_{j}/4})\prod_{\substack{p \in I_{j}}}\Big (1+\frac {k}{p}+ O(\frac 1{p^2}) \Big ) .
\end{align*}

  We then conclude that
\begin{align*}
\begin{split}
 & \sum_{\substack{n=\prod_jn_j }}\frac {\log n}{n}\prod^{\mathcal{J}}_{j=1}\Big(1-  2^{\Omega(n_j)-\ell_j/2}\Big ) \frac {k^{\Omega(n_j)}}{g(n_j)} b_j(n_j) \prod_{\substack{p \in I_j \\ p | n_j}}\Big (1 +O(\frac 1p) \Big ) \\
\geq &  \prod^{\mathcal{J}}_{\substack{j=1}}\Big ( 1-2^{-\ell_{j}/4} \Big ) \prod_{p \in \bigcup^{\mathcal{J}}_{j=1}I_j}\Big ( 1+\frac {k}{p}+ O(\frac 1{p^2}) \Big ) \Big (\sum_{p \in \bigcup^{\mathcal{J}}_{j=1}I_j}\frac {k\log p}{p}+O(\frac 1{p^2}) \Big ).
\end{split}
\end{align*}

   We apply the above estimation  into \eqref{loglowerbound} and apply \eqref{Sexp}, \eqref{SRexp} together with Lemma \ref{RS} to conclude that
\begin{align*}
  S_1 \gg & T(\log T)^{k^2+1} .
\end{align*}
   This completes the proof of the proposition.

\subsection{Proof of Proposition \ref{Prop6}}
\label{sec: proof of Prop 6}

   We denote the left side expression in \eqref{L2estmation} by $S_2$ and apply Cauchy’s integral formula for derivatives to see that
\begin{align*}
   S_2 &  = \int^T_1|\zeta'(\half+it)|^{2} |\mathcal{N}(\half+it, k-1)|^2 dt =
\int^T_1 \Big|\frac {1}{2\pi i}\int_{\mathcal{C}_1}\frac {\zeta(\half+\alpha+it)}{\alpha^2}d\alpha \Big |^{2} |\mathcal{N}(\half+it, k-1)|^2 dt,
\end{align*}
where $\mathcal{C}_1$ denotes the positively oriented circle in the complex plane centered at $0$ of radius $R=(\log T)^{-1}$. We then apply the Cauchy-Schwarz inequality to the integral over $\alpha$ above to deduce that
\begin{align}
\label{S2bound}
\begin{split}
   S_2 \leq & (\frac 1{2\pi})^2
\int^T_1 \Big|\int_{\mathcal{C}_1}\frac {1}{\alpha^4}d\alpha \Big |\Big |\int_{\mathcal{C}_1}|\zeta(\half+\alpha+it)|^2d\alpha\Big | \Big|\mathcal{N}(\half+it, k-1) \Big|^2 dt  \\
\leq & (\frac 1{2\pi})^2 R^{-2} \max_{|\alpha|=R}\int^T_1|\zeta(\half+\alpha+it)|^2|\mathcal{N}(\half+it, k-1)|^2 dt.
\end{split}
\end{align}

  We denote the last integral above by $I$ and we fix an $\alpha=\beta+i\gamma$ with $\beta, \gamma \in \mc$ such that $|\alpha|=R$ to estimate it. Without loss of generality, we may assume that $\beta \leq 0$ and we apply Cauchy’s residue theorem to deduce that
\begin{align*}
   I &  = \frac{1}{2 \pi i} \int_{\mathcal{C}_2} |\zeta(s)|^2  |\mathcal{N}(s-\alpha, k-1)|^2 \, ds,
\end{align*}
where $\mathcal{C}_2$ consists of line segments from $\half+\beta+(1+\gamma)i$ to $\half+(1+\gamma)i$, then from $\half+(1+\gamma)i$ to $\half+(T+\gamma)i$ and lastly from $\half+(T+\gamma)i$ to $\half+\beta+(T+\gamma)i$.

The integration on the on the horizontal edges of the contour can be estimated to be $O(T^{1-\varepsilon})$ using \eqref{Nnormbound} and the convexity bound for $\zeta(s)$  (see \cite[Exercise 3, p.
  100]{iwakow}) that asserts
\begin{align*}
\begin{split}
  \zeta(s) \ll & \left( 1+|s| \right)^{\frac {1-\Re(s)}{2}+\varepsilon}, \quad 0 \leq \Re(s) \leq 1,
\end{split}
\end{align*}

   We then deduce that
\begin{align*}
   I = I_{R}+ O(T^{1-\varepsilon}),
\end{align*}
 where
\begin{align*}
\label{IRint}
  I_{R} = \frac{1}{2 \pi i} \int_{1+\gamma}^{T+\gamma} |\zeta(\half+it)|^2|\mathcal{N}(\half+it-\alpha, k-1)|^2 dt.
\end{align*}

  We now apply arguments similar to the proof of \cite[Proposition 2]{H&Sound} to deduce that for $T$ large enough,
\begin{align*}
  I_{R} \ll T(\log T)^{k^2}.
\end{align*}

   We apply the above estimation in \eqref{S2bound} to conclude that
\begin{align*}
  S_2 \ll & T(\log T)^{k^2+2} .
\end{align*}
   This completes the proof of the proposition.

\vspace*{.5cm}

\noindent{\bf Acknowledgments.}  P. G. is supported in part by NSFC grant 11871082.

\bibliography{biblio}
\bibliographystyle{amsxport}

\vspace*{.5cm}

\end{document}